\documentclass[graybox]{svmult}

\usepackage{mathptmx} 
\usepackage{helvet}   
\usepackage{courier}  
\usepackage{graphicx} 
\usepackage{siunitx}
\usepackage{amsbsy}
\usepackage{amsmath,bm}
\usepackage{amssymb}
\usepackage{amsfonts}
\usepackage{multirow}
\usepackage[bottom]{footmisc}
\usepackage{cite}
\usepackage{url}

\spnewtheorem{assumption}{Assumption}{\bfseries}{\itshape}

\begin{document}

\title*{Stability and convergence of the penalty formulation for nonlinear magnetostatics}
\author{Herbert Egger 
\and
Felix Engertsberger 
\and
Klaus Roppert 
}
\institute{Herbert Egger \at RICAM Linz, Austria, 
\email{herbert.egger@ricam.oeaw.ac.at}
  \and Felix Engertsberger \at NUMA, JKU Linz, Austria 
\email{felix.engertsberger@jku.at}
  \and Klaus Roppert \at IGTE, TU Graz, Austria,\email{klaus.roppert@tugraz.at}}

\maketitle
 
\abstract*{The magnetostatic field distribution in a nonlinear medium amounts to the unique minimizer of the magnetic coenergy over all fields that can be generated by the same current. 
This is a nonlinear saddlepoint problem whose numerical solution can in principle be achieved by mixed finite element methods and appropriate nonlinear solvers.
The saddlepoint structure, however, makes the solution cumbersome.
A remedy is to split the magnetic field into a known source field and the gradient of a scalar potential which is governed by a convex minimization problem. 
The penalty approach avoids the use of artificial potentials and Lagrange multipliers and leads to an unconstrained convex minimization problem involving a large parameter.  
We provide a rigorous justification of the penalty approach by deriving error estimates for the approximation due to penalization. 
We further highlight the close connections to the Lagrange-multiplier and scalar potential approach. 
The theoretical results are illustrated by numerical tests for a typical benchmark problem.}

\abstract{The magnetostatic field distribution in a nonlinear medium amounts to the unique minimizer of the magnetic coenergy over all fields that can be generated by the same current. 
This is a nonlinear saddlepoint problem whose numerical solution can in principle be achieved by mixed finite element methods and appropriate nonlinear solvers.
The saddlepoint structure, however, makes the solution cumbersome.
A remedy is to split the magnetic field into a known source field and the gradient of a scalar potential which is governed by a convex minimization problem. 
The penalty approach avoids the use of artificial potentials and Lagrange multipliers and leads to an unconstrained convex minimization problem involving a large parameter.  
We provide a rigorous justification of the penalty approach by deriving error estimates for the approximation due to penalization. 
We further highlight the close connections to the Lagrange-multiplier and scalar potential approach. 
The theoretical results are illustrated by numerical tests for a typical benchmark problem.
}

\section{Introduction}
\label{egger:sec:intro}

The static distribution of magnetic fields and fluxes in a nonlinear anhysteretic medium can be described by the equations of magnetostatics~\cite{egger:Meunier}
\begin{alignat}{2}
\operatorname{curl} \mathbf{h} &= \mathbf{j} \quad \text{in } \Omega,   \qquad &\mathbf{b} &= \partial_\mathbf{h} w_*(\mathbf{h}), \label{egger:eq:1}\\
\operatorname{div} \mathbf{b} &= 0 \quad \text{in } \Omega, \qquad & \mathbf{b} \cdot \mathbf{n} &= 0 \quad \text{on } \partial\Omega. \label{egger:eq:2}
\end{alignat}
Here $\mathbf{h}$ and $\mathbf{b}$ denote the magnetic field and flux intensities, $\mathbf{j}$ is the given current density, and $\mathbf{n}$ is the normal vector at the boundary $\partial \Omega$. 
Further $\partial_\mathbf{h} w_*(\mathbf{h})$ denotes the derivative of the magnetic coenergy density $w_*(\mathbf{h})$ which allows to describe nonlinear anhysteretic material behavior~\cite{egger:Silvester1991}. 
For inhomogeneous materials, the coenergy density may additionally depend on the spatial coordinate.

\bigskip 
\noindent 
\textbf{Equivalent minimization problem.}
It is well-known~\cite{egger:Kikuchi1989} that solutions of the system \eqref{egger:eq:1}--\eqref{egger:eq:2} can be characterized equivalently by the 
minimization problem
\begin{align} \label{egger:eq:3}
\min_{\mathbf{h} \in H(\operatorname{curl})} \int_\Omega w_*(\mathbf{h}) \, dx \qquad \text{s.t.} \quad \operatorname{curl} \mathbf{h} = \mathbf{j} \quad \text{in } \Omega.
\end{align}
A standard way to treat the constraints in this system is provided by the Lagrange formalism~\cite{egger:Boffi}. The corresponding optimality conditions amount to
\begin{alignat}{2}
\partial_\mathbf{h} w_*(\mathbf{h}) - \operatorname{curl} \mathbf{a} &= \mathbf{0} \quad &&\text{in } \Omega, \label{egger:eq:4}\\
\operatorname{curl} \mathbf{h} &= \mathbf{j} \quad &&\text{in } \Omega, \label{egger:eq:5}
\end{alignat}
complemented by appropriate boundary conditions. 
From \eqref{egger:eq:4} and \eqref{egger:eq:1}, one can deduce that $\mathbf{b}=\operatorname{curl} \mathbf{a}$, i.e., the Lagrange multiplier $\mathbf{a}$ has the meaning of a magnetic vector potential. 
Additional gauging conditions are thus required in three space dimensions to guarantee its uniqueness~\cite{egger:GiraultRaviart,egger:Monk}.
Together with the saddlepoint structure already present in \eqref{egger:eq:4}--\eqref{egger:eq:5}, this complicates the numerical solution.

\bigskip 
\noindent 
\textbf{Reduced scalar potential approach.}
An alternative and more convenient way to deal with the constraints in \eqref{egger:eq:3} is to split the admissible fields as 
\begin{align} \label{egger:eq:6}
\mathbf{h} = \mathbf{h}_s - \nabla \psi, 
\end{align}
with given source field $\mathbf{h}_s$, satisfying $\operatorname{curl} \mathbf{h}_s = \mathbf{j}$, and reduced scalar potential $\psi$. This allows to rewrite \eqref{egger:eq:3} into an equivalent unconstrained minimization problem 
\begin{align} \label{egger:eq:7}
\min_{\psi \in H^1(\Omega)/\mathbb{R}} \int_\Omega w_*(\mathbf{h}_s - \nabla \psi) \, dx.
\end{align}
Here $H^1(\Omega)/\mathbb{R} = \{v \in H^1(\Omega) : \int_\Omega v \, dx=0\}$ is the space of scalar potentials with zero average; the latter condition is required to ensure uniqueness of the potential.
The strong form of the optimality system for this problem reads 
\begin{align} \label{egger:eq:8}
\operatorname{div}(\partial_\mathbf{h} w_*(\mathbf{h}_s-\nabla \psi)) &= 0 \qquad \text{in } \Omega,
\end{align}
again complemented by appropriate boundary conditions.
Let us note that the system \eqref{egger:eq:8} can also be obtained directly from \eqref{egger:eq:1}--\eqref{egger:eq:2} by eliminating $\mathbf{b}$ and $\mathbf{h}$ using~\eqref{egger:eq:6}. This approach is also well suited for numerical approximation~\cite{egger:Biro93,egger:Engertsberger}.

\bigskip 
\noindent 
\textbf{Penalty approach.}
A third possible way to deal with the constraints in \eqref{egger:eq:3} is via penalization~\cite{egger:Bandelier1994}; also see \cite{egger:Meunier}. This leads to the unconstrained minimization problem
\begin{align} \label{egger:eq:9}
\min_{\mathbf{h} \in H(\operatorname{curl})} \int_\Omega w_*(\mathbf{h})  + \frac{1}{2\varepsilon} |\operatorname{curl} \mathbf{h} - \mathbf{j}|^2 \, dx.
\end{align}
The existence of a unique solution $\mathbf{h}^\varepsilon$ here can be guaranteed for any $\varepsilon>0$, but note that the minimizer $\mathbf{h}^\varepsilon$ satisfies Amp\'ere's law $\operatorname{curl}\mathbf{h} = \mathbf{j}$ only approximately.

\bigskip 
\noindent 
\textbf{Scope.}
In this paper, we provide a rigorous justification of the penalty approach for nonlinear magnetostatics by deriving quantitative estimates for the approximation error introduced through penalization. 
Related error estimates for abstract regularized linear saddlepoint problems can be found in~\cite{egger:GiraultRaviart}. In a similar spirit, time dependent linear eddy current problems were analyzed recently~\cite{egger:Bermudez2020}.
The error estimates derived in this manuscript extend these previous results to the nonlinear magnetostatic setting. The proofs of our main results further clarify the close connection of the penalty method to the Lagrange multiplier and scalar potential approaches. 

\bigskip 
\noindent 
\textbf{Outline.}
In Section~\ref{egger:sec:main}, we introduce our notation and basic assumptions, and then state the main result of the paper. 
Its proof is given in Section~\ref{egger:sec:proof}. 
The possible generalization of the results and proofs to the discrete setting are briefly discussed in Section~\ref{egger:sec:discussion}.
For illustration of our theoretical results, we present some numerical tests in Section~\ref{egger:sec:num}, which also demonstrate the feasibility and efficiency of the penalty approach in combination with higher order finite element approximations. 

\section{Assumptions and main result}
\label{egger:sec:main}
By $H^1(\Omega)$, $H(\operatorname{curl};\Omega)$, $H(\operatorname{div};\Omega)$, and $L^2(\Omega)$, we denote the usual function spaces arising in the context of electromagnetic field problems; see e.g.~\cite{egger:GiraultRaviart,egger:Monk}.
For brevity, the symbol $\Omega$ will be omitted in most cases. 
We write $\|\cdot\|_X$ for the norm of the space $X$ and denote by $\langle a, b\rangle_\Omega = \int_\Omega a \cdot b \, dx$ the scalar product of $L^2(\Omega)^d$, $d \ge 1$. 
The following assumptions on the problem data will be utilized throughout the manuscript. 
\begin{assumption} \label{egger:ass:1}
$\Omega \subset \mathbb{R}^d$, $d=2,3$ is a bounded Lipschitz domain and simply connected. 
The coenergy density $w_* : \Omega \times \mathbb{R}^d \to \mathbb{R}$ is piecewise smooth with respect to~$x$, and for all $x \in \Omega$ the function $w_*(x,\cdot) : \mathbb{R}^d \to \mathbb{R}$ is smooth, strongly convex, and has bounded derivatives, i.e., 
\begin{align}
\langle \partial_\mathbf{h} w_*(x,\mathbf{u}) - \partial_\mathbf{h} w_*(x,\mathbf{v}), \mathbf{u}-\mathbf{v} \rangle &\ge \gamma \, |\mathbf{u}-\mathbf{v}|^2 \label{egger:eq:ass1}\\
|\partial_\mathbf{h} w_*(x,\mathbf{u})-\partial_\mathbf{h} w_*(x,\mathbf{v})| &\le L  \, |\mathbf{u}-\mathbf{v}| \label{egger:eq:ass2}
\end{align}
for all $\mathbf{u},\mathbf{v}\in \mathbb{R}^d$ and $x \in \Omega$ with uniform constants $L$, $\gamma>0$.
Here $\langle \cdot,\cdot\rangle$ and $|\cdot|$ are the Euclidean scalar product and norm on $\mathbb{R}^d$. 
Finally, $\mathbf{j} \in L^2(\Omega)^d$ with $\operatorname{div}\mathbf{j}=0$. % and $\mathbf{j} \cdot \mathbf{n}=0$ on $\partial\Omega$. 
\end{assumption}
The conditions of this assumption cover rather general inhomogeneous, nonlinear, and anisotropic materials, as well as permanent magnets. With minor modifications of our arguments below, the conditions on the topology of the domain $\Omega$ and the regularity of the coenergy density $w_*(\cdot)$ could be further relaxed to some extent. 
The above assumptions already allow us to state and prove our main result. 
\begin{theorem} \label{egger:thm:1}
Let Assumption~\ref{egger:ass:1} hold. Then the variational problems~\eqref{egger:eq:3} and~\eqref{egger:eq:7} each have a unique solution, and
\begin{align}
\mathbf{h} = \mathbf{h}_s - \nabla \psi, \qquad \mathbf{b} = \partial_\mathbf{h} w_*(\mathbf{h})
\end{align}
amounts to the unique (weak) solution of \eqref{egger:eq:1}--\eqref{egger:eq:2}. 
For any $\varepsilon>0$, also problem \eqref{egger:eq:9} has a unique solution $\mathbf{h}^\varepsilon$ and $\|\mathbf{h}^\varepsilon\|_{H(\operatorname{curl})} \le C \|\mathbf{j}\|_{L^2}$.
With $\mathbf{b}^\varepsilon = \partial_\mathbf{h} w_*(\mathbf{h}^\varepsilon)$, we have
\begin{align} \label{egger:eq:13} 
\|\mathbf{h}^\varepsilon - \mathbf{h}\|_{H(\operatorname{curl})} + \|\mathbf{b}^\varepsilon - \mathbf{b}\|_{H(\operatorname{div})} \le C' \varepsilon \|\mathbf{j}\|_{L^2}.
\end{align}
Moreover, the constants $C$, $C'$ in these estimates do not depend on the parameter $\varepsilon$.
\end{theorem}
The error in the penalty approximation can thus be controlled by choosing the regularization parameter $\varepsilon>0$ small, respectively, the penalty parameter $\frac{1}{\varepsilon}$ in \eqref{egger:eq:9} large enough.
The dependence on $\varepsilon$ hence has to be made clear in all estimates required for the proof of Theorem~\ref{egger:thm:1}. 
The analysis will further clarify the detailed connection between the penalty formulation \eqref{egger:eq:7} and the reduced scalar potential approach \eqref{egger:eq:9}. 

\section{Proof of the main theorem}
\label{egger:sec:proof}
By Assumption~\ref{egger:ass:1}, 
the system \eqref{egger:eq:9} amounts to a convex optimization problem, and the first order optimality condition
\begin{align} \label{egger:eq:14}
\langle \partial_\mathbf{h} w_*(\mathbf{h}^\varepsilon), \mathbf{v}\rangle_\Omega + \tfrac{1}{\varepsilon} \langle \operatorname{curl} \mathbf{h}^\varepsilon, \operatorname{curl} \mathbf{v} \rangle_\Omega &= \tfrac{1}{\varepsilon}\langle \mathbf{j}, \operatorname{curl} \mathbf{v}\rangle_\Omega \qquad \forall \mathbf{v}\in H(\operatorname{curl})
\end{align}
completley characterizes the minimizers of \eqref{egger:eq:9}; see e.g.~\cite{egger:Boffi,egger:Zeidler2b}. 
Existence of a unique solution to this variational problem can be shown by Zarantonello's fixed-point theorem~\cite{egger:Zeidler2b}. 
The existence of a unique solution to  problems \eqref{egger:eq:3} and \eqref{egger:eq:7}, their equivalence, as well as a uniform bound $\|\mathbf{h}\|_{H(\operatorname{curl})} \le C \|\mathbf{j}\|_{L^2}$ can be shown by similar arguments; see \cite{egger:Engertsberger} for the details. 
The fact that both approaches yield the solution to the original problem \eqref{egger:eq:1}--\eqref{egger:eq:2} follows by elementary arguments.

\medskip 
\noindent 
\textbf{Auxiliary results and rescaling.}
For proving the two estimates of Theorem~\ref{egger:thm:1}, 
let us recall from \cite{egger:GiraultRaviart} that any function $\mathbf{v}\in H(\operatorname{curl})$ can be decomposed uniquely as 
\begin{align} \label{egger:eq:15}
\mathbf{v}= \mathbf{z}' - \nabla \psi' 
\quad \text{with} \quad \mathbf{z}' \perp \nabla H^1
\end{align}
and $\psi' \in H^1/\mathbb{R}$. 
Note that $\operatorname{curl} \mathbf{z}' = \operatorname{curl} \mathbf{v}$ and recall the Poincar\'e inequality
\begin{align} \label{egger:eq:16}
c_P \|\mathbf{z}'\|_{H(\operatorname{curl})} \le  \|\operatorname{curl} \mathbf{z}'\|_{L^2} \qquad \forall \mathbf{z} \perp \nabla H^1,
\end{align}
which follows from general principles of functional analysis~\cite[Thm.~4.6]{egger:Arnold}.
By further splitting also $\mathbf{h}^\varepsilon = \mathbf{z}^\varepsilon - \nabla \psi^\varepsilon$ in the same way, we can decompose \eqref{egger:eq:14} into a system of two equations and, similar to \cite{egger:Eller2017}, we can rescale the second of these identities by the parameter $\epsilon$. This leads to the equivalent set of equations
\begin{alignat}{2}
\langle \partial_{\mathbf{h}} w_*(\mathbf{z}^\varepsilon - \nabla \psi^\varepsilon),\nabla \psi' \rangle_\Omega &= 0  \quad &\forall \psi' \in H^1/\mathbb{R}, \label{egger:eq:17}\\
\varepsilon \langle \partial_{\mathbf{h}} w_*(\mathbf{z}^\varepsilon - \nabla \psi^\varepsilon),\mathbf{z}'\rangle_\Omega + \langle \operatorname{curl} \mathbf{z}^\varepsilon, \operatorname{curl} \mathbf{z}'\rangle_\Omega 
&= \langle \mathbf{j}, \operatorname{curl} \mathbf{z}'\rangle_\Omega \quad &\forall \mathbf{z}' \perp \nabla H^1. \label{egger:eq:18}
\end{alignat}
Note that, after this rescaling, we can now set $\varepsilon=0$, which leads to the limit system 
\begin{alignat}{2}
\langle \partial_{\mathbf{h}} w_*(\mathbf{z} - \nabla \psi),\nabla \psi' \rangle_\Omega &= 0  \quad &\forall \psi' \in H^1(\Omega)/\mathbb{R},\label{egger:eq:19}\\
\langle \operatorname{curl} \mathbf{z}, \operatorname{curl} \mathbf{z}'\rangle_\Omega &= \langle \mathbf{j}, \operatorname{curl} \mathbf{z}'\rangle_\Omega \quad &\forall \mathbf{z}' \perp \nabla H^1.\label{egger:eq:20}
\end{alignat}
A quick inspection of the equations reveals that $\operatorname{curl} \mathbf{z} = \mathbf{j}$ by \eqref{egger:eq:20}, hence $\mathbf{z}$ amounts to a specific choice of the source field $\mathbf{h}_s$ in \eqref{egger:eq:6}. 
Furthermore, equation~\eqref{egger:eq:19} can be interpreted as the weak form of equation~\eqref{egger:eq:8}. 
Hence $\mathbf{h} = \mathbf{z} - \nabla \psi$ and $\mathbf{b} = \partial_\mathbf{h} w_*(\mathbf{h})$ amount to the unique weak solution of the system~\eqref{egger:eq:1}--\eqref{egger:eq:2}. 

\medskip 
\noindent 
\textbf{Error estimate.}
By substracting the  identities \eqref{egger:eq:18} and \eqref{egger:eq:20}, one can verify that
\begin{align} \label{egger:eq:21}
\varepsilon \langle \partial_\mathbf{h} w_*(\mathbf{h}^\varepsilon) - \partial_\mathbf{h} w_*(\mathbf{h}) , \mathbf{z}'\rangle_\Omega + \langle \operatorname{curl}(\mathbf{z}^\varepsilon - \mathbf{z}), \operatorname{curl} \mathbf{z}'\rangle_\Omega &= -\varepsilon \langle \partial_\mathbf{h} w_*(\mathbf{h}),\mathbf{z}'\rangle_\Omega
\end{align}
for all $\mathbf{z}' \perp \nabla H^1$. 
Choosing $\mathbf{z}' = \mathbf{z}^\varepsilon - \mathbf{z}$, using the properties of $\partial_\mathbf{h} w_*(\cdot)$, a Poincar\'e inequality~\cite{egger:GiraultRaviart}, and the bound $\|\mathbf{h}\|_{H(\operatorname{curl})} \le C \|\mathbf{j}\|_{L^2}$ then leads to 
\begin{align} \label{egger:eq:22}
c_P \|\mathbf{z}^\varepsilon - \mathbf{z}\|_{H(\operatorname{curl})} \le \|\operatorname{curl} (\mathbf{z}^\varepsilon - \mathbf{z})\|_{L^2} \le \frac{\varepsilon}{c_P} \|\partial_\mathbf{h} w_*(\mathbf{h})\|_{L^2} \le C \varepsilon \|\mathbf{j}\|_{L^2}.
\end{align}
From the identities \eqref{egger:eq:17} and \eqref{egger:eq:19}, on the other hand, we get 
\begin{align} \label{egger:eq:23}
\langle \partial_\mathbf{h} w_*(\mathbf{z}^\varepsilon - \nabla \psi^\varepsilon) - \partial_\mathbf{h} w_*(\mathbf{z} - \nabla \psi), \nabla \psi'\rangle_\Omega 
= 0
\end{align}
for all test functions $\psi' \in H^1/\mathbb{R}$. 
Using the monotonicity of $\partial_\mathbf{h} w_*(\cdot)$, one then obtains $\|\nabla \psi^\varepsilon - \nabla \psi\|_{L^2(\Omega)} \le c'' \|\mathbf{z}^\varepsilon - \mathbf{z}\|_{L^2(\Omega)}$. Together with the previous inequalities, this already yields the estimate $\|\mathbf{h}^\varepsilon - \mathbf{h}\|_{H(\operatorname{curl})} \le C \varepsilon \|\mathbf{j}\|_{L^2}$. 
The uniform bound for $\mathbf{h}^\varepsilon$ then follows from that for $\mathbf{h}$ and the triangle inequality.
From $\mathbf{b}^\varepsilon=w_*(\mathbf{h})^\varepsilon$, $\mathbf{b}=w_*(\mathbf{h})$, and the Lipschitz continuity of $\partial_\mathbf{h} w_*$, we further obtain 
\begin{align} \label{egger:eq:24}
\|\mathbf{b}^\varepsilon - \mathbf{b}\|_{L^2} \le L \|\mathbf{h}^\varepsilon - \mathbf{h}\|_{L^2} \le C'' \varepsilon \|\mathbf{j}\|_{L^2}.
\end{align}
Using~\eqref{egger:eq:17} and \eqref{egger:eq:19}, one can finally see that $\operatorname{div} \mathbf{b}^\varepsilon=\operatorname{div} \mathbf{b}=0$, which yields the bound for the second term in \eqref{egger:eq:13}, and completes the proof.
\hfill $\square$

\section{Observations and further results}
\label{egger:sec:discussion}

The analysis presented in the previous section not only provides a rigorous justification of the penalty (regularized $\mathbf{h}$-field) approach, but it also highlights its close connection to the reduced scalar potential formulation; see \eqref{egger:eq:17}--\eqref{egger:eq:18} and \eqref{egger:eq:19}--\eqref{egger:eq:20}.
A similar connection can be made to the Lagrange multiplier approach. By adding a term $\varepsilon \mathbf{a}$ to \eqref{egger:eq:5}, we obtain a regularized saddlepoint system
\begin{align}
\partial_\mathbf{h} w_*(\mathbf{h}^\varepsilon) - \operatorname{curl} \mathbf{a}^\varepsilon &= 0 \quad \text{in } \Omega, \label{egger:eq:25}\\
\operatorname{curl} \mathbf{h}^\varepsilon + \varepsilon \mathbf{a}^\varepsilon &= \mathbf{j} \quad \text{in } \Omega. \label{egger:eq:26}
\end{align}
Expressing $\mathbf{a}^\varepsilon$ via \eqref{egger:eq:26} and inserting the result into \eqref{egger:eq:25} leads to 
\begin{align}
\partial_\mathbf{h} w_*(\mathbf{h}^\varepsilon) + \tfrac{1}{\varepsilon}\operatorname{curl} \mathbf{h}^\varepsilon &= \tfrac{1}{\varepsilon} \mathbf{j} \qquad \text{in } \Omega, 
\end{align}
which amounts to the strong form of \eqref{egger:eq:14}. 
A careful inspection of the main arguments further reveals that all results carry over almost verbatim also to $H(\operatorname{curl})$-conforming approximations of the variational problems, e.g., by finite elements or iso-geometric analysis.
The error estimate of Theorem~\ref{egger:thm:1} thus remains valid on the discrete level with constants $C$, $C'$ that are independent of the discretization parameters, e.g., the mesh size $h$ or the polynomial degree $p$. 
Some post-processing is required, of course, on the discrete level to obtain approximations for $\mathbf{b}^\varepsilon$ and $\mathbf{b}$ in $H(\operatorname{div})$.
Similar arguments can also be used to establish estimates for the discretization error for the regularized $\mathbf{h}$-field approach that are independent of the regularization parameter $\varepsilon$. 
A detailed analysis in these directions will be provided elsewhere. In the following section, we however illustrate these facts by some numerical tests.

\section{Numerical tests}
\label{egger:sec:num}

For illustration of the theoretical results, we now briefly present some numerical results obtained for a typical test problem involving linear and nonlinear media. 

\bigskip 
\noindent 
\textbf{Linear and nonlinear material laws.}
The magnetic response of air and copper is assumed linear with $\mathbf{b} = \mu_0 \mathbf{h}$ which can be expressed as
\begin{align} \label{egger:bh_air}
\mathbf{b} = \partial_{\mathbf{h}} w_*(\mathbf{h}) 
\qquad \text{with} \qquad 
w_*(\mathbf{h}) = \frac{\mu_0}{2} |\mathbf{h}|^2.
\end{align}
Since the material is oviously isotropic, we can express $w_*(\mathbf{h})=\widetilde w_*(|\mathbf{h})|$ with the scalar function $\widetilde w_* : \mathbb{R}_+ \to \mathbb{R}_+$ here defined by $\widetilde w_*(s) = \frac{\mu_0}{2} s^2$. 
In a similar manner, the nonlinear isotropic response of iron can be described as 
\begin{align} \label{egger:bh_iron}
\mathbf{b} 
= \partial_{\mathbf{h}} w_*(\mathbf{h}) 
= \widetilde w_*'(|\mathbf{h}|) 
\frac{\mathbf{h}}{|\mathbf{h}|}.
\end{align}
This amounts to $\mathbf{b} = \mu_{ch}(|\mathbf{h}|) \, \mathbf{h}$
with $\mu_{ch}(\mathbf{h}) = \widetilde w_*'(|\mathbf{h}|)/|\mathbf{h}| $ denoting the chord permeability. 
By taking the norm on both sides of \eqref{egger:bh_iron}, we see that $|\mathbf{b}| = \widetilde w_*'(|\mathbf{h}|)$, i.e., the derivative $\widetilde w_*'(|\mathbf{h}|)$ can be obtained directly from the usual B-H curves. 
For our numerical tests, we take the B-H data from the TEAM~13 benchmark problem~\cite{egger:team13}, which are interpolated piecewise linearly. The function $\widetilde w_*(|\mathbf{h}|)$ is then obtained by integration of this spline function. 
The conditions for our analysis are thus satisfied.

\bigskip 
\noindent 
\textbf{Implementation.}
For the numerical solution of \eqref{egger:eq:14}, we consider a Galerkin approximation with $H(\operatorname{curl})$-conforming finite elements of higher order. Curved simplicial elements are used for meshing  the geometry. The nonlinear systems are solved by a Newton-method with linesearch for which global convergence can be established~\cite{egger:Engertsberger}.
All computations are performed in Netgen/NGSolve~\cite{egger:ngsolve}. 

\begin{figure}[ht]
\centering
\hfill
\includegraphics[height=0.28\textwidth]{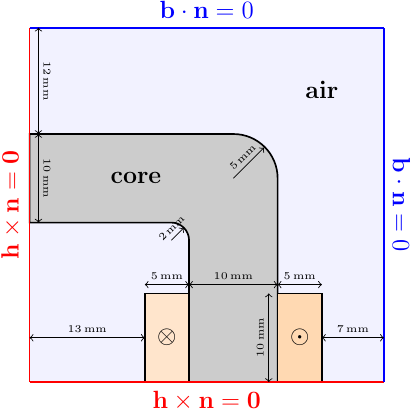}
\hfill
\includegraphics[trim=10 26 33 15,clip,height=0.28\textwidth]{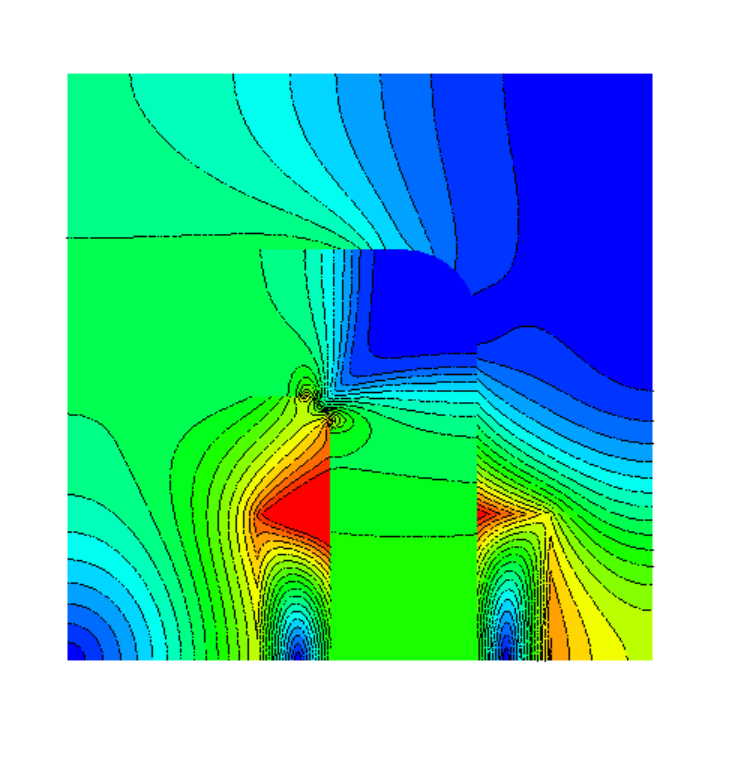}
\hfill
\includegraphics[height=0.28\textwidth, trim=33 31 10 18,clip]{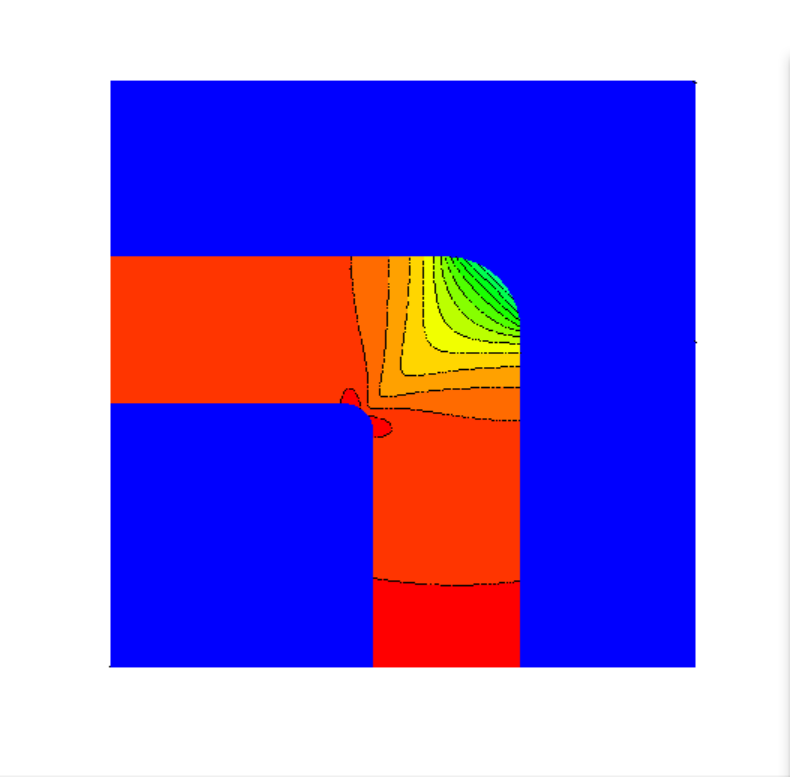}
\hfill
\caption{Left: Geometry of test problem. 
Middle: Magnetic field density $|\mathbf{h}|$ (values $0 \text{--} 30.000 \,\unit{\ampere}/\unit{\meter}$) computed by the finite element approxmation of \eqref{egger:eq:14} on a mesh with $\text{nt}=58.240$ triangles and with polynomial order $\text{p}=2$.
Right: Corresponding flux density $|\mathbf{b}| = \widetilde w_*'(|\mathbf{h})|$ (values $0$--$2 \, \unit{\tesla}$).  
\label{egger:fig:1}}
\end{figure}

\bigskip 
\noindent 
\textbf{Test problem.}
The geometry depicted in Figure~\ref{egger:fig:1} models a quarter of the two-dimensional cross-section of a transformer; see \cite{egger:friedrich2019} for details. 
Since we are in a two dimensional setting, $\mathbf{j}=(0,0,j_3)$ and $\mathbf{h}=(h_1,h_2,0)$ with $h_1$, $h_2$, $j_3$ independent of the $z$-coordinate.
We are interested in the magnetostatic fields generated by an excitation current with density $j_3=\pm 10^{7} \, \unit{\ampere\per\square\meter}$.
For small values of the regularization parameter $\epsilon>0$, the numerical approximations obtained for \eqref{egger:eq:14} are visually indistinguishable from those obtained by a corresponding vector potential formulation, which may serve as an informal validation of the results; also see Table~\ref{egger:tab:1}.

\begin{table}[ht!]
\centering
\setlength{\tabcolsep}{0.4em}
\def\arraystretch{1.3}
\footnotesize
\begin{tabular}{c||c|c|c|c|c||c}
$\epsilon_0$ & $10^{-1}$ & $10^{-2}$ & $10^{-3}$ & $10^{-4}$ & $10^{-5}$ & disc.\,err.\\
\hline 
$\frac{\|\mathbf{h}_h - \mathbf{h}_h^\epsilon\|_{L^2}}{\|\mathbf{h}_h\|_{L^2}}$ & $6.41 \cdot 10^{-1}$ & $7.27 \cdot 10^{-2}$ & $7.35 \cdot 10^{-3}$ & $7.35 \cdot 10^{-4}$ & $7.31 \cdot 10^{-5}$ &  $1.42 \cdot 10^{-4}$
\\
\hline
$\frac{\|\mathbf{b}_h - \mathbf{b}_h^\epsilon\|_{L^2}}{\|\mathbf{b}_h\|_{L^2}}$ & $1.68 \cdot 10^{-1}$ & $1.47 \cdot 10^{-2}$ & $1.45 \cdot 10^{-3}$ & $1.48 \cdot 10^{-4}$ & $1.55 \cdot 10^{-5}$ & $1.78 \cdot 10^{-4}$ 
\end{tabular}
\medskip
\caption{Convergence of regularized approximations $\mathbf{h}_h^\epsilon$, $\mathbf{b}_h^\epsilon$ for \eqref{egger:eq:14} to reference solution $\mathbf{h}_h$, $\mathbf{b}_h$ of system \eqref{egger:eq:6}--\eqref{egger:eq:7} with regularization parameter $\epsilon = \epsilon_0 \ell^2/\mu_0$, where $\ell=40 \, \unit{mm}$ is the domain size. 
%
%The last column contains estimates for the relative discretization errors. 
\label{egger:tab:1}}
\end{table}

\noindent 
\textbf{Computational results.}
To illustrate the estimates of Theorem~\ref{egger:thm:1}, we compare in Table~\ref{egger:tab:1} the numerical solution of the regularized H-field approximation \eqref{egger:eq:14} with regularization parameter $\epsilon>0$ to the corresponding numerical approximation of the reduced scalar potential formulation \eqref{egger:eq:6}--\eqref{egger:eq:7}, which amounts to the limit $\epsilon=0$ in~\eqref{egger:eq:14}. 
As mentioned in Section~\ref{egger:sec:discussion}, the estimates of  Theorem~\ref{egger:thm:1} carry over almost verbatim to the discrete level. 
In the last column, we also display the relative differences to the numerical solution obtained by the vector potential ($\mathbf{a}$-field) formulation on the same mesh, which may serve as an estimate for the discretization errors. 

\bigskip 
\noindent 
\textbf{Discussion.}
The results of our computations clearly demonstrate linear convergence of the approximation error with respect to the regularization parameter $\epsilon$, which is in perfect agreement with the theoretical predictions made in Theorem~\ref{egger:thm:1}.
Already for a moderate size of the penalization parameter $\epsilon$, the discretization errors, which due to the use of a higher order approximation, curved elements, and adapted meshes, are very small here, start to dominate the effect of the errors introduced by penalization. 
Similar results were also obtained for test problems in three dimensions, where the use of parameter robust iterative solvers is required for an efficient realization of the penalty approach. Investigations in this direction are left for future research.

\begin{acknowledgement}
The authors are grateful for financial support by the international FWF/DFG funded Collaborative Research Centre CREATOR (TRR361/SFB-F90).
\end{acknowledgement}

\end{document}